\documentclass[11pt]{article}%
\usepackage{fullpage}
\usepackage{amsmath}
\usepackage{amsfonts}
\usepackage{amssymb}
\usepackage{graphicx}
\usepackage{hyperref}%
\providecommand{\U}[1]{\protect\rule{.1in}{.1in}}
\begin{document}

\begin{center}
{\Large \textbf{Simple finite elements and multigrid for efficient
mass-consistent wind downscaling in a~coupled fire-atmosphere model}}
\newline\ \newline{J. Mandel\,${}^{1}$, A. Farguell\,${}^{2}$, A. K.
Kochanski\,${}^{2}$, D. V. Mallia\,${}^{3}$, K. Hilburn\,${}^{4}$}\newline%
{${}^1\,$University of Colorado Denver, Denver, CO\\
${}^2$San Jos\'e State University, San Jos\'e, CA\\
${}^3$University of Utah, Salt Lake City, UT\\
${}^4$Colorado State University, Fort Collins, CO}%

\end{center}

\section{Introduction}

In the coupled atmosphere-fire model
WRF-SFIRE~\cite{Mandel-2014-NFW,Mandel-2011-CAF}, the Weather Research
Forecasting (WRF) model~\cite{Skamarock-2008-DAR} runs at 300m--1km horizontal
resolution, while the fire model runs at the resolution of 30m or finer. The
wind has a fundamental effect on fire behavior and the topography details have
a strong effect on the wind, but WRF\ does not see the topography on the fire
grid scale. We want to downscale the wind from WRF\ to account for the
fine-scale terrain.\ For this purpose, we fit the wind from WRF with a
divergence-free flow over the detailed terrain. Such methods, called
mass-consistent approximations, were originally proposed on regular
grids~\cite{Sherman-1978-MMW,Singh-2008-EQF} for urban and complex terrain
modeling, with terrain and surface features modeled by excluding entire grid
cells from the domain. For fire applications,
WindNinja~\cite{Wagenbrenner-2016-DSW} uses finite elements on
a~terrain-following grid. The resulting equations are generally solved by
iterative methods such as SOR, which converge slowly, so use of GPUs is of
interest~\cite{Bozorgmehr-2020-CIF}. A~multigrid method with
a~terrain-following grid by a~change of coordinates was proposed
in~\cite{Wang-2005-AMM}.

The method proposed here is to be used in every time step of WRF-SFIRE in the
place of interpolation to the fire model grid. Therefore, it needs to have the
potential to (1) scale to hundreds or thousands of processors using WRF
parallel infrastructure~\cite{Wang-2019-AUG}; (2) scale to domains size at
least 100km by 100km horizontally, with $3000\times3000\times15$ grid cells
and more; (3) have reasonable memory requirements per grid point; (4) not add
to the cost of the time step significantly when started from the solution in
the previous time step; and, (5) adapt to the problem automatically, with
minimum or no parameters to be set by the user.

\section{Finite element formulation}

\label{sec:formulation}

Given vector field $\boldsymbol{u}_{0}$ on domain $\Omega\subset\mathbb{R}%
^{d}$, subset $\Gamma\subset\partial\Omega$, and $d\times d$ symmetric
positive definite coefficient matrix $\boldsymbol{A}=\boldsymbol{A}\left(
\boldsymbol{x}\right)  $, we want to find the closest divergence-free vector
field $\boldsymbol{u}$ by solving the problem%
\begin{equation}
\min_{\boldsymbol{u}}\frac{1}{2}\int\limits_{\Omega}\left(  \boldsymbol{u}%
-\boldsymbol{u}_{0}\right)  \cdot\boldsymbol{A}\left(  \boldsymbol{u}%
-\boldsymbol{u}_{0}\right)  d\boldsymbol{x}\text{\quad subject to
}\operatorname{div}\boldsymbol{u}=0 \text{ in } \Omega\text{ and
}\boldsymbol{u}\cdot\boldsymbol{n}=0\text{ on }\Gamma,
\label{eq:continuous-min-constr}%
\end{equation}
where $\Gamma$ is the bottom of the domain (the surface), and $\boldsymbol{A}%
\left(  \boldsymbol{x}\right)  $ is a $3\times3$ diagonal matrix with penalty
constants $a_{1}^{2},a_{2}^{2},a_{3}^{2}$ on the diagonal. Enforcing the
constraints in (\ref{eq:continuous-min-constr}) by a~Lagrange multiplier
$\lambda$, we obtain the solution $\left(  \boldsymbol{u},\lambda\right)  $ as
a stationary point of the Lagrangean%
\begin{equation}
\mathcal{L}\left(  \boldsymbol{u},\lambda\right)  =\frac{1}{2}\int%
\limits_{\Omega}\boldsymbol{A}\left(  \boldsymbol{u}-\boldsymbol{u}%
_{0}\right)  \cdot\left(  \boldsymbol{u}-\boldsymbol{u}_{0}\right)
d\boldsymbol{x}+\int\limits_{\Omega}\lambda\operatorname{div}\boldsymbol{u}%
d\boldsymbol{x}-\int\limits_{\Gamma}\lambda\boldsymbol{n}\cdot\boldsymbol{u}%
d\boldsymbol{s}. \label{eq:continuous-lagrangean}%
\end{equation}
Eliminating $\boldsymbol{u}$ from the stationarity conditions $\partial
\mathcal{L}(\boldsymbol{u},\lambda)/\partial\lambda=0$ and $\partial
\mathcal{L}(\boldsymbol{u},\lambda)/\partial\boldsymbol{u}=0$ by
\begin{equation}
\boldsymbol{u}=\boldsymbol{u}_{0}+\boldsymbol{A}^{-1}\operatorname{grad}%
\lambda\label{eq:recover-u}%
\end{equation}
leads to the generalized Poisson equation for Lagrange multiplier $\lambda$,%
\begin{equation}
-\operatorname{div}\boldsymbol{A}^{-1}\operatorname{grad}\lambda
=\operatorname{div}\boldsymbol{u}_{0}\text{ on }\Omega,\quad\lambda=0\text{ on
}\partial\Omega\setminus\Gamma,\text{ \quad}\boldsymbol{n\cdot A}%
^{-1}\operatorname{grad}\lambda=-\boldsymbol{n\cdot u}_{0}\text{ on }\Gamma.
\label{eq:laplace}%
\end{equation}
Multiplication of (\ref{eq:laplace}) by a test function $\mu$, $\mu=0$ on
$\partial\Omega\setminus\Gamma$, and integration by parts yields the
variational form to find $\lambda$ such that $\lambda=0$ on $\partial
\Omega\setminus\Gamma$ and
\begin{equation}
\int_{\Omega}\boldsymbol{A}^{-1}\operatorname{grad}\lambda\cdot
\operatorname{grad}\mu\,d\boldsymbol{x}=-\int_{\Omega}\operatorname{grad}%
\mu\cdot\boldsymbol{u}_{0}d\boldsymbol{x} \label{eq:variational}%
\end{equation}
for all $\mu$ such that $\mu=0$ on $\partial\Omega\setminus\Gamma$. The
solution is then recovered from (\ref{eq:recover-u}). We proceed formally
here; see~\cite{Juarez-2012-MCW} for a~different derivation
of~(\ref{eq:variational}) in a~functional spaces setting.

The variational problem (\ref{eq:variational}) is discretized by standard
isoparametric 8-node hexahedral finite elements, e.g.,~\cite{Hughes-1987-FEM}.
The integral on the left-hand side of (\ref{eq:variational}) is evaluated by
tensor-product Gauss quadrature with two nodes in each dimension, while for
the right-hand side, one-node quadrature at the center of the element is
sufficient. The same code for the derivatives of a finite element function is
used to evaluate $\operatorname{grad}$ $\lambda$ in (\ref{eq:recover-u}) at
the center of each element.

The unknown $\lambda$ is represented by its values at element vertices, and
the wind vector is represented naturally by its values at element centers. No
numerical differentiation of $\lambda$ from its nodal values, computation of
the divergence of the initial wind field $\boldsymbol{u}_{0}$, or explicit
implementation of the boundary condition on $\operatorname{grad}\lambda$ in
(\ref{eq:laplace}) is needed. These are all taken care of by the finite
elements naturally.

\section{Multigrid iterations}

The finite element method for (\ref{eq:variational}) results in a system of
linear equations $Ku=f$. The values of the solution are defined on a~grid,
which we will call a~\emph{fine grid}. One cycle of the multigrid method
consists of several iterations of a basic iterative method, such as
Gauss-Seidel, called a~\emph{smoother}, followed by a \emph{coarse-grid
correction}. A prolongation matrix $P$ is constructed to interpolate values
from a coarse grid, in the simplest case consisting of every other node, to
the fine grid. For a given approximate solution $u$ after the smoothing, we
seek an improved solution in the form $u+Pu_{c}$ variationally, by solving%
\begin{equation}
P^{\top}K\left(  u+Pu_{c}\right)  =P^{\top}f \label{eq:correction-system}%
\end{equation}
for $u_{c}$, and obtain the coarse-grid correction procedure as%
\begin{align}
f_{c}=P^{\top}\left(  f-Ku\right)  \qquad &  \text{form the coarse right-hand
side}\nonumber\\
K_{c}=P^{\top}KP\qquad &  \text{form the coarse stiffness matrix}\nonumber\\
K_{c}u_{c}=f_{c}\qquad &  \text{solve the coarse-grid problem}%
\label{eq:coarse-problem}\\
u\leftarrow u+Pu_{c}\qquad &  \text{insert the coarse-grid correction}%
\nonumber
\end{align}
The coarse grid correction is followed by several more smoothing steps, which
completes the multigrid cycle.

In the simplest case, $P$ is a linear interpolation and the coarse stiffness
matrix $K_{c}$ is the stiffness matrix for a coarse finite element
discretization on a grid with each coarse-grid element taking the place of a
$2\times2\times2$ agglomeration of fine-grid elements. That makes it possible
to apply the same method to the coarse-grid problem (\ref{eq:coarse-problem})
recursively. This process creates a hierarchy of coarser grids. Eventually,
the coarsest grid problem is solved by a direct method, or one can just do
some more iterations on it.

Multigrid methods gain their efficiency from the fact that simple iterative
methods like Gauss-Seidel change the values of the solution at a node from
differences of the values between this and neighboring nodes. When the error
values at neighboring nodes become close, the error can be well approximated
in the range of the prolongation $P$ and the coarse-grid correction can find
$u_{c}$ such that $u+Pu_{c}$ is a much better approximation of the solution.
For analysis of variational multigrid methods and further references, see
\cite{Bank-1981-OOP,Mandel-1987-VMT}.

Multigrid methods are very efficient. For simple elliptic problems, such as
the Poisson equation on a regular grid, convergence rates of about $0.1$
(reduction of the error by a factor of $10$) at the cost of $4$ to $5$
Gauss-Seidel sweeps on the finest grid are expected \cite{Brandt-1977-MLA}.
However, the convergence rates get worse on more realistic grids, and
adaptations are needed. We choose as the smoother vertical sweeps of
Gauss-Seidel from the bottom up to the top, with red-black ordering
horizontally into $4$ groups. For the base method, we use $2\times2\times2$
coarsening and construct $P$ so that the vertices of every $2\times2\times2$
agglomeration of elements interpolate to the fine-grid nodes in the
agglomeration, with the same weights as the trilinear interpolation on
a~regular grid. The interpolation is still trilinear on a~stretched grid, but
only approximately trilinear on a deformed terrain-following grid.

The base method works as expected as long as some grid layers are not tightly
coupled. If they are, we mitigate the slower convergence by semicoarsening
\cite{Morano-1998-CSU}: After smoothing, the error is smoother in the tightly
coupled direction(s), which indicates that we should not coarsen the other
direction(s). When the grid is stretched vertically away from the ground, the
nodes are relatively closer and thus tightly coupled in the horizontal
direction. Similarly, when the penalty coefficient $a_{3}$ in the vertical
direction is larger than $a_{1}$ and $a_{2}$ in the horizontal directions, the
neighboring nodes in the vertical direction are tightly coupled numerically.
The algorithm to decide on coarsening we use is: Suppose that the penalty
coefficients are $a_{1}=a_{2}=1$ and $a_{3}\geq1$, and at the bottom of the
grid, the grid spacing is $h_{1}=h_{2}$ (horizontal) and $h_{3}$ (vertical).
If $h_{3}/(h_{1}a_{3})>1/3$, coarsen in the horizontal directions by $2$,
otherwise do not coarsen. Then, replace $h_{1}$ and $h_{2}$ by their new
values, corsened (multiplied by 2)\ or not, and for every horizontal layer
from the ground up, if $h_{3}/(h_{1}a_{3})]<3$, coarsen about that layer
vertically, otherwise do not coarsen. This algorithm retains the coarse grids
as logically cartesian, which is important for computational efficiency and
keeping the code simple, and it controls the convergence rate to remain up to
about $0.28$ with four smoothing steps per cycle.

\section{Conclusion}

We have presented a simple and efficient finite element formulation of
mass-consistent approximation, and a multigrid iterative method with adaptive
semicoarsening, which maintains the convergence of iteration over a range of
grids and penalty coefficients. A prototype code is available at \url{https://github.com/openwfm/wrf-fire-matlab/tree/femwind/femwind}.

\textbf{Acknowledgement:} This work has been supported by NSF grant
ICER-1664175 and NASA grant 80NSSC19K1091.


\end{document}